
\input amstex
\documentstyle{amsppt}

\nologo
\topspace{2.50in}
\document
\centerline{ON REMOVABLE SETS FOR SOBOLEV SPACES IN THE PLANE}

\vskip .25in
\centerline{Peter W. Jones\footnote{Supported by NSF Grant
DMS-8916968}}
\centerline{Department of Mathematics}
\centerline{Yale University}
\centerline{New Haven, CT 06520}

\newpage

\baselineskip 20pt
Let $K$ be a compact subset of $\bar{\bold C} ={\bold R}^2$ and let $K^c$ denote
its complement.  We say $K\in HR$, $K$ is holomorphically removable,
if whenever $F:\bar{\bold C}\to\bar{\bold C}$ is a homeomorphism and
$F$ is holomorphic off $K$, then $F$ is a M\"obius transformation.
By composing with a M\"obius transform, we may assume $F(\infty
)=\infty$.  The contribution of this paper is to show that a large
class of sets are $HR$.  Our motivation for these results is that
these sets occur naturally (e.g. as certain Julia sets) in dynamical
systems, and the property of being $HR$ plays an important role in
the Douady-Hubbard description of their structure.  (See \cite{4}.)

To prove that the sets in question are $HR$ we establish what may be
a stronger result.  A compact set $K$ is said to be removable for
$W^{1,2}$ if every $f$ which is continuous on ${\bold R}^2$ and in the
Sobolev space $W^{1,2}(K^c)$ (one derivative in $L^2$ on $K^c$) is
also in $W^{1,2}(\bold R^2)$.  It is a fact that if $K$ is removable
for $W^{1,2}$, $K$ is $HR$.  We do not know the answer to the
following question:

If $K$ is $HR$, is $K$ removable for $W^{1,2}$?

To prove the fact we first show that the two dimensional Lebesgue
measure of $K$, $\vert K\vert$, is zero.  If not let
$F_n=\frac{1}{\pi z}*(e^{in(x+y)}\chi_K(z))$.  then $\underset
n\rightarrow\infty\to\lim\Vert F_n\Vert_{L^{\infty}(\bold R^{2})}=0$
and $F_n$ is continuous.  Since $\vert\bar\partial F_n\vert =\chi_n$,
$\Vert\bar\partial F_n\Vert_{L^{2}(\bold R^{2})}=\vert K\vert^{1/2}$.
 On the other hand, $F_n(z)\to 0$ for $z\not\in K$, so $L^4$ bounds
on convolution with $\frac{1}{\pi z^{2}}$ when combined with
H\"older's inequality show $\Vert F'_n\Vert_{L^{2}(K^{c})}\to 0$.
(See \cite{11}, Chapter 1.)  Taking a sum of functions like the
$F_n$, we obtain a globally continuous $F\in W^{1,2}(K^c)$, $F\not\in
W^{1,2}(\bold R^2)$.  Now using the fact that $\vert K\vert =0$, we
deduce $K\in HR$.  Take a homeomorphism $F$ with $F'(\infty )=1$.
Then $f(z)=F(z)-z\in
W^{1,2}(K^c)$ because integrating $\vert F'\vert^2$ gives the area of
the image.  Now $f\in W^{1,2}(\bold R^2)$ and $\bar\partial f=0$
except on a set of measure zero implies (Weyl's lemma) $f$ is
holomorphic.  Therefore $F(z)=z+a$.

We recall some elementary facts concerning $HR$.  If $\vert K\vert >0$,
it follows from
the ``measurable Riemann mapping theorem'' (see \cite{1}) that there
is a nontrivial quasiconformal mapping $F$ which is holomorphic off
$K$.  (Thus $K\notin HR$.)  If $K$ has Hausdorff dimension less than
1, $\text{Dim}(K)<1$, the fact that $K\in HR$ follows from the Cauchy
integral formula (Painleve's theorem).  Similarly, if $K$ is a
rectifiable curve, Morera's theorem implies $K\in HR$.  Kaufman
\cite{7} has produced examples of curves where $\text{Dim}(K)=1$ but
$K\notin HR$.  The ``difficult'' case is the one that occurs in
conformal dynamics:  $K$ is connected and has some ``fractal''
properties.  (The case of ``pure'' Cantor-type sets is easy; they are
$HR$.  By a pure Cantor set, we mean e.g. one arising from a Cantor
construction with a constant ratio of dissection, or the Julia set for
$z^2+c$ where $c$ is not in the Mandelbrot set.)  We also point out
that the case where $K$ is a quasicircle
seems to be folklore - again, $K\in HR$.  That the property of being
$HR$ is related to quasiconformal mappings is seen from the following

\proclaim{Remark}  $K$ is $HR$ if and only if whenever $F$ is a
homeomorphism of $\bar{\bold C}$ which is $M$ quasiconformal on
$K^c$, $F$ is globally quasiconformal (and hence $M$ quasiconformal).
 (See \cite{8}, page 200.)
\endproclaim

To prove the remark, first assume that $K$ is $HR$.  By the
measurable Riemann mapping theorem there is a globally quasiconformal
mapping $G$ such that $G\circ F$ is holomorphic off $K$.  Since
$G\circ F$ is a M\"obius transformation and $\vert K\vert =0$, $F$ is
globally $(M)$ quasiconformal.  For the other direction, standard
$L^p$ estimates (see \cite{1}) show that necessarily $\vert K\vert
=0$.  If $F$ is a homeomorphism which is analytic off $K$, $F$ is
globally quasiconformal and hence $(\vert K\vert =0)$ a M\"obius
transformation.

For $\Omega$ be a domain on the Riemann sphere and let $z_0\in\Omega$.  Then $\Omega$ is a John domain
(with center $z_0$) if there is $\varepsilon >0$ such that for all $z_1\in
\Omega$ there is an arc $\gamma\subset\Omega$ which connects $z_0$ to
$z_1$ and has the property that
$$
d(z)\geq\varepsilon d(z,z_1),\quad z\in\gamma.
$$
Here $d(z,z_1)$ is the chordal distance from $z$ to $z_1$ and $d(z)$
is the chordal distance of $z$ to $\partial\Omega$  We call such an
arc $\gamma$ a John arc.  In this paper we
will choose coordinates so that $z_0=\infty$, and this allows us to
replace $d(z),d(z,z_1)$ by the corresponding Euclidean distances.
The property of being a John domain is preserved under globally
quasiconformal mappings.  If $\Omega$ is a simply connected John
domain, it is easy to show that the arc $\gamma$ may be taken to be
the hyperbolic geodesic from $z_0$ to $\infty$.  (See \cite{9} for an
exposition of properties of John domains.)  The main result of this
paper is

\proclaim{Theorem 1}  If $\Omega$ is a John domain and
$K=\partial\Omega$, then $K$ is removable for $W^{1,2}$.
\endproclaim

Notice that the hypothesis demands that $K=\partial\Omega$, but says
{\it nothing} about the other components of $\bold C\backslash K$.
This is because the hypothesis will be seen to force some geometry on
those other components.  (For example, the interior of a cardioid is
a John domain while the exterior is not.  The parabolic basin for
$z^2+\frac{1}{4}$ is also a John domain, while the basin for $\infty$
- the exterior domain - is not.)  It is of some philosophical interest to
note the similarities between Theorem 1 and the results of \cite{6} on
extension problems for Sobolev spaces.

Since the John condition is quasiconformally invariant, we obtain
directly (see also ``Remark'')

\proclaim{Corollary 1}  If $\Omega$ is a John domain and
$K=\partial\Omega$, any global homeomorphism which is quasiconformal
off $K$ is globally quasiconformal (with the same constant of
quasiconformality).
\endproclaim

We say that a polynomial $P(z)$ is subhyperbolic on its Julia set $J$
if there is a metric $\lambda (z)\vert dz\vert$ such that $\lambda
(z) -\underset j\to\sum \vert z-z_j\vert^{-\alpha_{j}}$ is $C^\infty$
for some numbers $\alpha_j <1$, and $P(z)$ is hyperbolic on $J$ in
the metric $\lambda$.  In other words there are $c,\varepsilon >0$
such that for all $n\geq 1$,
$$
\lambda (z)^{-1}\lambda (P_n(z))\vert\frac{d}{dz} P_n(z)\vert\geq
c(1+\varepsilon )^n.
$$
Here $P_n(z)=P\circ\cdots\circ P(z)$ is the $n^{\text{th}}$ iterate
of $P$.  (This definition may be a bit restrictive, but it is all we
will need for this paper.)  The following question is open:
$$
\text{If}~ J~\text{is the Julia set for a polynomial, is}~J\in HR?
$$
It is proven in \cite{3} that whenever a polynomial $P(z)$ is
subhyperbolic on its Julia set $J$, then $A_\infty$, the basin of
attraction at $\infty$ for $P$, is a John domain.  Since $J=\partial
A_\infty$, we obtain

\proclaim{Corollary 2}  If $P(z)$ is subhyperbolic on its Julia
set $J$, then $J\in HR$.
\endproclaim

The corollary answers a question of A. Douady and J. Hubbard and was
the starting point of this investigation.  
Douady posed the question to the author for the particular
(subhyperbolic) case where
$P(z)=z^2+c$ has the (Misiurewicz) property that the origin is
preperiodic but not
periodic (e.g. $z^2+i$).  This case is not fundamentally different
for the general case of subhyperbolic polynomials.  An amusing
feature of our proof is that the Julia set for a Misiurewicz point
(from the family $z^2+c$) is actually
easier to deal with than those arising from the hyperbolic case.
(When $K^c=\Omega$, our argument is a bit simpler.  The arguments of
Sections 5 and 6 are not needed.)

The proof of Theorem 1 starts by proving it in the case where
$\Omega$ is simply connected on $\bar{\bold C}$, i.e. $K$ is
connected.  The general case then follows from

\proclaim{Theorem 2}  If $\Omega$ is an $(\varepsilon )$ John domain,
there is a $(c(\varepsilon ))$ John domain $\Omega '$ with $\Omega '$
simply connected and
$$
\partial\Omega\subset\partial\Omega '.
$$
\endproclaim

While the proof of Theorem 2 is perhaps not immediately obvious, it
turns out to follow from a simple construction with planar graphs.

Section 2 contains background material, and Sections 3-7 are devoted
to the proof of Theorem 1.  The idea is to redefine $F$ near $K$ so
that it is $C^\infty$ near $K$ and so that the Sobolev norm does not
change much.  Theorem 2 is proven in Section 8.

\vskip .25in

\centerline{\bf {\S2.  Background Material}}

Let $F\in W^{1,2}(K^c)$ be continuous on $\bold C$.  An easy argument
with the Dirichlet principle shows that to prove $F\in W^{1,2}(\bold
R^2)$ it is sufficient to treat the case we now assume, where $F$ is
harmonic near $K$.  We also assume the reader is familiar with
elementary properties of logarithmic capacity, which we denote by
$\text{Cap}(\circ )$.  See e.g. \cite{10} for the first two of the
next three lemmata.
Let $f:\bold D\to\bold C$ be univalent, $f(0)=0, f'(0)=1$.  Then $f$
has a Fatou extension to ${\bold T} =\partial{\bold D}$ and this
extension is always defined except on a set of capacity (and hence
Lebesgue measure) zero.  In our applications, all image domains will
have locally connected boundaries, and hence $f$ will be continuous
on $\bar{\bold D}$.  The following 
results are due to Beurling.  (See e.g. \cite{10} for Lemmata 2.1 and
2.2.)   The values of $c$ below are various
universal constants.

\proclaim{Lemma 2.1}  If $E\subset\bold T$,
$$
\text{Cap}(f(E))\geq c~\text{Cap}(E)^2.
$$
\endproclaim

\proclaim{Lemma 2.2}  Let $g_\theta =f(\{re^{i\theta} :0\leq r<1\})$.
Then if $\ell (\cdot )$ denotes arclength,
$$
\text{Cap}(\{ e^{i\theta}:\ell (g_\theta )>\lambda\})\leq
c\lambda^{-1/2}
$$
\endproclaim

\proclaim{Lemma 2.3}  Suppose $H$ is harmonic and continuous in $\bold
D$, and $(\vert\nabla H\vert^2=\vert H_x\vert^2+\vert H_y\vert^2)$,
$$
\iint_{\bold D}\vert\nabla H\vert^2dxdy=1.
$$
Then
$$
\text{Cap}(\{e^{i\theta}:\vert
H(e^{i\theta})-H(0)\vert\geq\lambda\})\leq ce^{-\pi\lambda^{2}}.
$$
\endproclaim

This last lemma can be found on page 30 of \cite{2}.  
We next require some elementary geometric facts about simply
connected John domains.  For the next result see \cite{5}.

\proclaim{Lemma 2.4}  If $g$ is a Poincar\'e geodesic from $\infty$
to $z_0\in\partial\Omega$ where $\Omega$ is an $(\varepsilon )$ John
domain, then $g$ is an arc of a $K(\varepsilon )$ quasicircle.
\endproclaim

Suppose now $\Omega$ is a bounded $(\varepsilon )$ John domain and
suppose the John center $z_0$ satisfies $d(z_0)=1$, where
$$
d(z)=~\text{distance}(z,\partial\Omega ).
$$
Then $\text{diameter}(\Omega )\sim 1$.  Let $f:\bold D\to\Omega ,f(0)=z_0$ be any choice of Riemann mapping, and define for $E\subset\partial\Omega$,
$$
\text{Cap}(E,z_0,\Omega
)\equiv~\text{Cap}(\{e^{i\theta}:f(e^{i\theta} )\in E\}).
$$

\proclaim{Lemma 2.5}  For any Borel set $E\subset\partial\Omega$,
$$
\text{Cap}(E,z_0,\Omega )\sim~\text{Cap}(E).
$$
\endproclaim

In the last line we mean that $A\sim B$ if there is a constant
$M=M(\varepsilon )$ such that
$$
M^{-1}A^M\leq B\leq MA^{1/M}.
$$
To prove the lemma let $G(z)=G(z,z_0)$ be Green's function for
$\Omega$ with pole at $z_0$.  Then it follows from the John condition
and the Koebe $\frac{1}{4}$ theorem that
$$
G(z)\geq cd(z)^\alpha~,~\alpha =\alpha (\varepsilon ),
$$
whenever $\vert z-z_0\vert\geq\frac{1}{2}$.  Suppose now that $z_j\in
E, z_j=f(\zeta_j),j=1,2$.  Fix a point $\zeta_3\in\bold D$ such that
$$
(1-\vert\zeta_3\vert
)\sim\vert\zeta_1-\zeta_3\vert\sim\vert\zeta_2-\zeta_3\vert\sim\vert\zeta_1-\zeta_2\vert
$$
and let $z_3=f(\zeta_3)$.  Then by the John condition
$$
\vert z_1-z_2\vert\leq Cd(z_3),
$$
while by our last estimate,
$$
d(z_3)\leq CG(z_3)^{1/\alpha}\sim C\vert\zeta_1-\zeta_2\vert^{1/\alpha}.
$$
In other words,
$$
\vert\zeta_{_{1}}-\zeta_{_{2}}\vert\geq c\vert z_1-z_2\vert^\alpha,
$$
and it follows from the definition of logarithmic capacity that
$$
\text{Cap}(E,z_0,\Omega )\geq M^{-1}\text{Cap}(E)^M.
$$
The other direction of the lemma follows from Lemma 2.1.\qed

\proclaim{Lemma 2.6}  Suppose $\Omega_j$ are $(\varepsilon )$ John
domains with centers $z_j, j=1,2$, and suppose $d(z_1),d(z_2)\sim 1$.
 Suppose also that $F$ is harmonic on $\Omega_1\cup\Omega_2$ and
continuous on $\bar\Omega_1\cup\bar\Omega_2$.  Then if
$E\subset\partial\Omega_1\cap\partial\Omega_2$ satisfies
$$
\text{Cap}(E)\geq\delta >0,
$$
there are geodesics $g_j\subset\Omega_j$ from $z_j$ to
$\partial\Omega_j$ such that $g_1$ and $g_2$ terminate at the same
point $\zeta\in\partial\Omega_1\cap\partial\Omega_2$, and
$$
\vert F(\zeta )-F(z_j)\vert\leq A(\varepsilon ,\delta
)(\iint_{\Omega_{j}}\vert\nabla F\vert^2dxdy)^{1/2},~j=1,2.
$$
\endproclaim

\demo{Proof}  Let
$$
E_j=\{z\in\partial\Omega_j:\vert
F(z)-F(z_j)\vert\geq\lambda(\iint_{\Omega_{j}}\vert\nabla
F\vert^2dxdy)^{1/2}\}.
$$
If $\lambda$ is large enough, Lemmata 2.3 and 2.5 show
$\text{Cap}(E_1\cup E_2)<\delta$.  Then $E\backslash (E_1\cup
E_2)\neq\phi$, so we may select $\zeta$ from that set.\qed

\vskip .25in

\centerline{{\bf \S3.  Quasicircles}}

We now give a quick outline of our proof for the case where $K$ is a
quasicircle.  This represents the only idea of the paper.  The rest
of the sections contain only technical arguments which make the same
philosophy work for the general case.

Let $\Omega_+$ and $\Omega_-$ denote respectively the unbounded and
bounded components of $\bar{\bold C}\backslash K$.  Fix two points
$z_\pm\in\Omega_\pm$ satisfying
$$
\delta (z_+)\sim\delta (z_-)\sim\vert z_+-z_-\vert\sim\delta,
$$
and build domains ${\Cal D}_\pm\subset\Omega_\pm$ which are bounded
by quasicircles and such that $\partial{\Cal D}_+\cap\partial{\Cal
D}_-$ is a subarc of $K$ with diameter $\sim\delta$.  The points
$z_\pm$ are made to be the ``centers'' of ${\Cal D}_\pm$.  Then by
Lemma 2.6 there is $A$ such that
$$
\vert F(z_+)-F(z_-)\vert\leq A(\iint_{{\Cal D}_{+}\cup{\Cal
D}_{-}}\vert\nabla F\vert^2dxdy)^{1/2}.
$$
Standard smoothing techniques now show there is $\tilde F\in
W^{1,2}(\bold R^2)$ such that $\tilde F =F$ outside of $K_\delta =\{
z:d(z)\leq\delta\}$, $\tilde F$ is $C^\infty$ near $K$, and
$$
\iint_{K_{\delta}}\vert\nabla\tilde F\vert^2dxdy\leq
c\iint_{K_{c\delta}}\vert\nabla F\vert^2dxdy.
$$
Sending $\delta$ to zero we see that $F\in W^{1,2}(\bold R^2)$ and
$$
\iint_{\bold R^{2}}\vert\nabla F\vert^2dxdy =\iint_{K^{c}}\vert\nabla
F\vert^2dxdy.
$$

If $F$ is $M$ quasiconformal on $K^c$, Lemma 2.2 and an argument
similar to the one above show that $F$ is globally quasiconformal.
The point of this vague remark is that, \it whatever argument we use\rm ,
it should show that $F$ being $M$ quasiconformal on $K^c$ implies $F$
is globally quasiconformal.  (See the ``Remark'' in Section 1.)

\vskip .24in
\centerline{{\bf \S4.  Some Geometry}}

In this section we construct certain domains related to a point
$x_0\in K$ and a scale $r$.  Since the John condition is scale
invariant, we may assume $x_0=0$ and $r=1$.  We will add to $K$
certain curves to obtain a new set $\hat K$ so that, in a certain
sense, $\bold C\backslash\hat K$ looks like a union of quasidisks of
diameter about $1$ (near $K$).

Let $f:\bold D^*=\{\vert z\vert >1\}\to\Omega$ be univalent with
$f(\infty )=\infty$.  Since $K=\partial\Omega$ is locally connected,
$f$ is continuous up to $\bold T$.  Select angles
$0=\theta_0<\theta_1<\theta_2<\cdots <\theta_N=2\pi$ so that
$$
\vert f(e^{i\theta})-f(e^{i\theta_{j}})\vert\leq
1~,~\theta_j\leq\theta\leq\theta_{j+1},
$$
and
$$
\vert f(e^{i\theta_{j}})-f(e^{i\theta_{j+1}})\vert\geq\frac{1}{2}.
$$
Now fix $M\geq 1$ and let $r_j<1$ be the largest value of $r$ so that
$$
\vert f(re^{i\theta_{j}})-f(e^{i\theta_{j}})\vert =M.
$$
Setting $L_j=\{re^{i\theta_{j}},r_j\leq r<1\}$ we see that
$$
\text{distance}(L_j,L_k)\geq c(1-r_j)~,~ j\neq k,\tag4.1
$$
for otherwise the John condition would be violated for the
corresponding geodesics in $\Omega$.

\proclaim{Lemma 4.1}  $\vert 1-r_j\vert\sim\vert
1-r_{j+1}\vert\sim~\text{distance}(L_j,L_{j+1})$.
\endproclaim

\demo{Proof}  We show that $\vert\theta_{j+1}-\theta_j\vert\leq
C(1-r_j)$.  The proof that $\vert\theta_{j+1}-\theta_j\vert\leq
C(1-r_{j+1})$ is the same.  The lemma will then follow from (4.1).
Let $I=\{e^{i\theta}:\theta_j\leq\theta\leq\theta_j+\pi\}$.  By
symmetry
$$
\omega (\zeta_j,I,\bold D^*)=\frac{1}{2},
$$
where $\zeta_j\equiv r_je^{i\theta_{j}}$. Here $\omega (z,E,{\Cal
D})$ denotes the harmonic measure at $z$ of $E\subset\partial{\Cal
D}$ in ${\Cal D}$.   By Beurling's so-called
$\frac{1}{2}$ theorem \cite{10}, if we set
$I_j=\{e^{i\theta}:\theta_j\leq\theta\leq\theta_{j+1}\}$,
$$
\omega (\zeta_j,I_j,\bold D )=\omega (f(\zeta_j),f(I_j),\Omega )\leq
CM^{-\frac{1}{2}},
$$
because diameter $(f(I_j))\leq 1$ and distance
$(f(\zeta_j),f(I_j))\geq M-1$.  Thus
$$
\omega (\zeta_j,I\backslash I_j,\bold D^*)\geq\frac{1}{4}
$$
if $M$ is large enough, and the lemma follows from simple estimates
on the Poisson kernel.\qed

Let ${\Cal D}_j$ be the domain bounded by $\bold T ,L_j,L_{j+1}$, and
the line segment $[\zeta_j,\zeta_{j+1}]$ and let
$\tilde\zeta_j=R_je^{i\varphi_{j}}$ where
$R_j-1=\frac{1}{2}~\text{min}(r_j-1,r_{j+1}-1)$, and
$\varphi_j=\frac{1}{2}(\theta_j+\theta_{j+1})$.  Then since we are
assuming diameter $(K)>>1$, each ${\Cal D}_j$ looks like a
quadrilateral (in $\bold D^*$ with one side on $\bold T$) with
bounded geometry.

\proclaim{Lemma 4.2}  $\Omega_j=f({\Cal D}_j)$ is an $(\varepsilon
')$ John domain with John center $z_j=f(\tilde\zeta_j)$.
\endproclaim

\demo{Proof}  Let $\zeta\in{\Cal D}_j$ and let $L=[\zeta
,\tilde\zeta_j]$ be the line segment from $\zeta$ to $\tilde\zeta_j$.
 Then $L\subset{\Cal D}_j$ and if $\zeta '\in L$, distance$(\zeta '
,\partial{\Cal D}_j)\geq c\vert\zeta '-\zeta\vert$.  (This follows
from the elementary geometry of ${\Cal D}_j$.)  Now if $L'$ is the
geodesic from $\zeta$ to $\infty$ in $\bold D^*$, $L'=\{R\zeta :R\geq
1\}$, $\rho (\zeta ',L')\leq C$ for all $\zeta '\in L$, where $\rho$
is the hyperbolic metric on $\bold D^*$.  The lemma now follows from
the John property on the arc $f(L')$ and the distortion theorem for
$f$.  The details are left to the reader.\qed

Lemma 4.2 is actually a special case of the following fact:

If $f:\bold D\to\Omega ,f(0)=z_0$, and $\Omega$ is an $(\varepsilon
)$ John domain with John center $z_0$, and if ${\Cal D}\subset\bold
D$ is a $(\delta )$ John domain with John center the origin, then
$f({\Cal D})$ is an $\eta (\varepsilon ,\delta )$ John domain with
John center $z_0$.

We leave a proof of this statement as an exercise for the reader.

At this point we remark that $\hat\Omega_j=$ interior of
$\bar\Omega_j$ is a $\delta (\varepsilon )$ quasicircle if
$\bar{\bold C}\backslash K=\Omega$.  (This is e.g. the case for the
Julia set corresponding to $z^2+i$.)  A most unfortunate complication
is that this statement is easily seen to be false if $\bar{\bold
C}\backslash K$ is allowed to have bounded components.  This
necessitates the technical construction of our next section.  The
reader interested only in the case where $\Omega =\bar{\bold
C}\backslash K$ may skip to Section 7, noticing that Proposition 6.1
has already been proven for quasicircles.

\vskip .25in
\centerline{\bf \S5.  Some Additional Curves}

We now add some additional curves to $K$.  Let ${\Cal O}_j$ be a
bounded component of $\bar{\bold C}\backslash K$.  Then by the
definition of the domains $\Omega_k$, each $\partial\Omega_k$
intersects $\partial{\Cal O}_j$ in either a connected set or the empty
set.  Let us for  the moment reorder the $\Omega_k$ so that
$\Omega_1,\dots ,\Omega_N$ are exactly those domains such
that $\partial\Omega_n\cap\partial{\Cal O}_j$ consists of more than one
point (and hence an arc).  Let $\delta >0$ be a small constant to be
fixed later and fix a Riemann mapping $f_j:\bold D\to{\Cal O}_j$ so
that $I_1,\dots ,I_N$ are intervals with
$f_j(I_n)=\partial\Omega_n\cap\partial{\Cal O}_j$.  By selecting
$f_j(0)$ to lie very close to $\partial\Omega_1\cap\partial{\Cal O}_j$
we may assume $\ell (I_1)\approx 2\pi$.  Let $T_1$ be the tent shaped
region bounded by $\bold T\backslash I_1$ and two straight lines in
$\bold D$ which intersect $\bold T\backslash I_1$ at angle $\delta$.
The $T_1$ is a ``thin sliver''.  Define ${\Cal U}_1=\hat{\Cal
U}_1=\bold D\backslash T_1$ so that $\partial{\Cal U}_1\cap\bold T
=I_1$.  For $n\geq 2$ let $L^1_n,L^2_n$ be the two lines which start
at the endpoints of $I_n$, go into $\bold D$, and make angle
$=\delta$ with $\bold T\backslash I_n$.  Let
$J_n=\{(1-\delta^{-1}\ell (I_n))e^{i\theta}:0\leq\theta\leq 2\pi\}$,
and let $\hat{\Cal U}_n$ be the domain bounded by the four arcs $I_n,
L^1_n,L^2_n,J_n$.  Then $\hat{\Cal U}_n$ almost fills up a rectangle
with length (along $\bold T$) $=\delta^{-2}\ell (I_n)$ and width (in
the direction orthogonal to $\bold T$) $=\delta^{-1}\ell (I_n)$.
Then by elementary estimates on the Poisson kernel,
$$
\{\omega (z,I_n,\bold D )\geq c_1\delta\}\subset\hat{\Cal
U}_n\subset\{\omega (z,I_n,\bold D )>c_2\delta\}.\tag5.1
$$

By reordering we may now assume that 
$$
\ell (I_2) \geq\ell (I_3)\geq\cdots\geq\ell (I_N).
$$
Define ${\Cal U}_n=\hat{\Cal U}_n\backslash\underset k=1\to{\overset
n\to\bigcup}{\Cal U}_k$ so that $\underset n=1\to{\overset
N\to\bigcup}{\Cal U}_n=\underset n=1\to{\overset
N\to\bigcup}\hat{\Cal U}_n$.  Recall that a domain ${\Cal U}$ is
called an $M$ Lipschitz domain if there is $z_0\in{\Cal U}$ and $R>0$
such that
$$
\partial{\Cal U}=\{z_0+Rr(\theta )e^{i\theta}:0\leq\theta\leq 2\pi\}
$$
where
$$
(1+M)^{-1}\leq r(\theta)\leq 1\quad\text{for all}~\theta
$$
and
$$
\vert r(\theta )-r(\theta ')\vert\leq M\vert\theta -\theta '\vert.
$$

\proclaim{Lemma 5.1}  ${\Cal U}_n$ is a $M(\delta )$ Lipschitz
domain, $1\leq n\leq N$.  Furthermore, if $\zeta\in\bold
D\cap\partial{\Cal U}_n$, there is $\varphi\in[0,2\pi ]$ such that
the line segment $\bar{\bold D}\cap\{\zeta +re^{i\theta}:r\geq 0\}$
lies in $\bar{\Cal U}_n$ and has endpoint on $I_n$ whenever
$$
\vert\theta -\varphi\vert\leq c\delta^2.
$$
\endproclaim

The proof of the lemma is an exercise in elementary geometry.  Now
let ${\Cal O}^n_j=f_j({\Cal U}_n)$.  If we consider any $\Omega_k$, we
have for each ${\Cal O}_j$, such that
$\partial{\Cal O}_j\cap\partial\Omega_k$ is an arc, obtained a domain
${\Cal O}_j^k\subset{\Cal O}_j$ (sometimes ${\Cal O}_j^k={\Cal O}_j$) with
the property that
$\partial{\Cal O}_j^k\cap\partial{\Cal O}_j\subset\partial\Omega_k$.  Let
${\Cal F}_k=\{{\Cal O}_j^k:\partial {\Cal O}_j 
\cap\partial\Omega_k~\text{is an arc}\}$ and let 
$\tilde\Omega_k=$ interior of closure of 
$\Omega_k\cup\underset{\Cal F}_k\to\bigcup{\Cal O}_j^k$.

\proclaim{Lemma 5.2}  $\partial\tilde\Omega_k$ is an $\eta
(\varepsilon ,\delta )$ quasicircle.
\endproclaim

\demo{Proof}  Let $\gamma_j^k={\Cal O}_j\cap\partial{\Cal O}_j^k$.  Then
$\partial\tilde\Omega_k\subset\partial\Omega_k\cup\underset{\Cal
F}_k\to\bigcup\gamma_j^k$.  We first claim $\tilde\Omega_k$ is an
$\eta (\varepsilon ,\delta)$ John domain.  It is only necessary to
find for every $z_0\in\partial\tilde\Omega_k$ an arc
$\gamma\subset\tilde\Omega_k$ which has endpoints $z_0$ and $z_k$
such that
$$
\text{distance}(z,\partial\tilde\Omega_k)\geq\eta\vert
z-z_0\vert~,~z\in\gamma.
$$
If $z_0\in\partial\Omega_k$ this is clear by Lemma 4.2.  We therefore
assume $z_0\in\gamma_j^k$ for some $j$.  By Lemmata 2.2, 2.3 and 5.1
there are angles $\varphi_{-1}<\varphi_0<\varphi_1$ such that
$\vert\varphi_\ell -\varphi_m\vert\sim\delta^2~,~\ell\neq m$,
such that
$$
\bold D\cap\{f_j^{-1}(z_0)+re^{i\theta},r>0\}\subset{\Cal U}_j,
$$
whenever $\varphi_{-1}\leq\theta\leq\varphi_1$, and such that
$$
\ell (\Gamma_m)\equiv\ell
(f_j(\{f_j^{-1}(z_0)+re^{i\varphi_{m}}:r>0\})\leq Cd(z_0).
$$
Furthermore, Lemma 5.1 allows us to assume that $\Gamma_m$ is a
Jordan arc and if $z\in\Gamma_\ell$ and $\vert
z-z_0\vert\geq\frac{1}{2} d(z_0)$,
$$
\text{distance}(z,\Gamma_m)\geq cd(z_0)~,~\ell\neq m.\tag5.2
$$
Now let $\delta_m$ be the endpoint of $\Gamma_m$ on
$\partial\Omega_k$ and let $\gamma_{-1}$ (resp. $\gamma_1$) be the
John geodesic from $\zeta_{-1}$ (resp. $\zeta_1$) in $\Omega_k$ to
$z_k$ (the John center of $\Omega_k$).  Then the curve $\gamma
=\Gamma_{-1}\cup\Gamma_1\cup\gamma_{-1}\cup\gamma_1$ surrounds
$\zeta_0$ and by the John condition on $\Omega_k$,
$$
\text{distance}(\zeta_0,\gamma )\geq cd(z_0).\tag5.3
$$
(Notice here that we are implicitly using the fact that $z\in{\Cal
O}_j^k$ implies $d(z)\leq C$.  This in turn follows from (5.1) and
either Lemma 2.2 or 2.3.)  Notice also that the interior of $\gamma$
must lie entirely in $\tilde\Omega_k$.

Let $\gamma_0$ be the John geodesic in $\Omega_k$ from $\zeta_0$ to
$z_k$.  We claim that the John condition for $\tilde\Omega_k$ holds
on $\Gamma_0\cup\gamma_0$.  First suppose that $z\in\Gamma_0$ and
$\vert z-z_0\vert\leq\frac{1}{2} d(z_0)$.  Then by the distortion
theorem for $f_j$,
$$
\text{distance}(z,\partial\tilde\Omega_k)\geq~\text{distance}(z,\gamma
)\geq c\vert z-z_0\vert.
$$
Now by inequality (5.2),
$$
d(z,\partial\tilde\Omega_k)\geq d(z,\gamma )\geq c\vert z-z_0\vert
$$
whenever $z\in\Gamma_0$ and $\vert z-z_0\vert\geq\frac{1}{2}d(z_0)$.
(Here we have used the John property on $\Omega_k$ to obtain distance
$(z,\gamma_{-1}\cup\gamma_1)\geq c\vert z-z_0\vert$.)

We must finally check the John condition on $\gamma_0$.  If
$z\in\gamma_0$ and $\vert z-\zeta_0\vert\leq cd(z_0)$, the inequality
for 
distance$(z,\partial\tilde\Omega_k)$ follows from (5.3) and the fact
that $\vert z_0-\zeta_0\vert\leq Cd(z_0)$.  If $z\in\gamma_0$ and
$\vert z-\zeta_0\vert >cd(z_0)$, the inequality for
distance$(z,\partial\tilde\Omega_k)$ follows from the John condition
distance$(z,\partial\Omega_k)\geq c\vert z-\zeta_0\vert$.  We have
thus established that $\tilde\Omega_k$ is a John domain.

We now claim that $G_k=\bar{\bold
C}\backslash\overline{(\tilde\Omega_k)}$ is a John domain.  We note
that by the definition of $\tilde\Omega_k$,
$\partial\tilde\Omega_k=\partial G_k$.  Now fix a point $z_0\in G_k$.

\enddemo
\enddemo
\enddemo
\enddemo
\subheading{Case A}

$z_0\in\overline{(\tilde\Omega_j)}$
for some
$j$.  (Then $j\neq k$).  First draw the John geodesic in $\Omega_j$
from $z_0$ to $z_j$.  We then draw the geodesic (in the Poincar\'e
metric of $\Omega$) from $z_j$ to $\infty$.  By the construction of
the domains $\Omega_j$ (Lemma 3.2) this is a John geodesic.  The
union of these two geodesics provides the arc joining $z_0$ to
$\infty$.

\subheading{Case B}  $z_0\in\Omega\backslash\underset
j\to\bigcup\bar\Omega_j$.  Let $\gamma$ be the Poincar\'e geodesic in
$\Omega$ from $z_0$ to $\infty$.  Then by Lemma 3.2, $\gamma$ is a
John geodesic in $G_k$.

\subheading{Case C}  $z_0\notin\Omega\cup\underset
j\to\bigcup\tilde\Omega_j$.  Then $z_0\in{\Cal O}_{j_{0}}$ for some
$j_0$.  Let
$$
A=\underset j\to\sup~\text{diameter}~\Omega_j,
$$
so that $A\sim 1$.  If $d(z_0)\geq 2A$ there is a half line $\gamma$
(to $\infty$ from $z_0$) which is a John geodesic in $G_k$.  If
$d(z_0)<2A$ there is a hyperbolic geodesic (which is also a John
arc in $G_k$) $\gamma_1$ from $z_0$ to $z_1\in{\Cal
O}_{j_{0},\ell}$ where $\ell\neq k, d(z)\geq 1$ on $\gamma_1$, and
$\ell (\gamma )\leq C$.  (This follows from the definition of the
domains ${\Cal O}_{j,k}$.)  By Case A there is a John geodesic
$\gamma_2$ from $z_1$ to $\infty$ in $G_k$.  The curve $\gamma
=\gamma_1\cup\gamma_2$ is the required John arc.

The proof is now completed by first observing that a simply connected
domain ${\Cal D}$ with $\partial{\Cal D}$ locally connected, ${\Cal
D}=~\text{Interior}(\bar{\Cal D})$, and $\bold C\backslash\bar{\Cal
D}$ connected is bounded by a Jordan curve, and then invoking the
following fact (see \cite{9}):

A Jordan domain ${\Cal D}$ is bounded by a quasicircle if and only if
${\Cal D}$ and $\bar{\bold C}\backslash\bar{\Cal D}$ are John domains.

\vskip .25in

\centerline{\bf \S6.  An Estimate on Capacity}

We now seek to imitate the proof given in Section 3.  What is
required is an estimate implying that $\vert F(z)-F(z_j)\vert$ is not
too large on $\partial\tilde\Omega_j$, except for a set of small
capacity.  While $\tilde\Omega_j$ is a quasidisk, $\tilde\Omega_j\cap
K^c$ is not necessarily connected.  This means we cannot simply apply
Lemma 2.3.  We state our result as a proposition; its proof will be
broken into several steps.  The result we state is far from optimal,
but it is all we need.  Let $\Tilde{\Tilde\Omega}_k$ be the domain
obtained by adding to $\tilde\Omega_k$ the set
$$
\underset j\to\bigcup \{z\in{\Cal O}_j:\rho (z,\partial{\Cal
O}_{j,k})<1\},
$$
where $\rho$ is the hyperbolic metric on ${\Cal O}_j$.  The domains
$\Tilde{\Tilde\Omega}_k$ then satisfy
$$
\sum\chi_{\Tilde{\Tilde\Omega}_{k}}\leq C.
$$

\proclaim{Proposition 6.1}  Suppose $H$ is continuous on the closure
of $\Tilde{\Tilde\Omega}_k$ and harmonic on
$\Tilde{\Tilde\Omega}_k\backslash K$.  Then if
$$
\underset\Tilde{\Tilde\Omega}_k\backslash K\to\iint \vert\nabla
H\vert^2dxdy = 1,
$$
we have the estimate
$$
\text{Cap}(\{ z\in\partial\tilde\Omega_k:\vert H(z)-H(z_k)\vert
>\lambda\})=o(1)
$$
as $\lambda\to\infty$.
\endproclaim

\demo{Proof}  Let
$E_1=\{z\in\partial\tilde\Omega_k\cap\partial\Omega_k:\vert
H(z)-H(z_k)\vert >\lambda\}$ and let
$E_2=\{z\in\partial\tilde\Omega_k\backslash\partial\Omega_k:\vert
H(z)-H(z_k)\vert >\lambda\}$.  Then by Lemmata 2.3 and 2.5,
$\text{Cap}(E_1)=o(1)$ as $\lambda\to\infty$, so it is sufficient to
show $\text{Cap}(E_2)=o(1)$ as $\lambda\to\infty$.

\enddemo
\subheading{Step 1.  Construction of Some Special Points}

Let $\{z_n\}$ be a collection of points in
$\partial\tilde\Omega_k\backslash\partial\Omega_k$ satisfying
$$
\vert z_n-z_m\vert\geq\frac{1}{4} d(z_n)~,~\forall n,m
$$
and
$$
\underset n\to\inf\vert z-z_n\vert\leq\frac{1}{2} d(z)~,~\forall
z\in\partial\tilde\Omega_k\backslash\partial\Omega_k.
$$

We will now form for each $z_n$ a point $z^*_n\in\Omega_k$.  For an
arbitrary point $z\in\partial{\Cal
O}_{j,k}\backslash\partial\Omega_k$ we let
$K_z=\{\zeta\in\partial\Omega_k:\vert z-\zeta\vert\leq 2d(z)\}$ so
that $\text{Cap}(K_z)\geq cd(z)$.  By the John condition and Lemma
2.5, there is a point $z^*\in\Omega_k$ such that $d(z)\sim
d(z^*)\sim\vert z-z^*\vert$ and there is a set $\tilde K_z\subset
K_z$ such that
$$
\text{Cap}(\tilde K_z,z^*,\Omega ),\text{Cap}(\tilde K_z,z,{\Cal
O}_j)\geq c.\tag6.1
$$

Denote by $f$ a Riemann mapping from $\Omega_k$ to $\bold D$ with
$f(z_k)=0$, where $z_k$ is the ``center'' of $\Omega_k$.  We can move
$z^*$ slightly so that $f(z^*)$ has the form
$$
f(z^*)=(1-2^{-\ell})\text{exp}\{im2^{-\ell}\pi\}\tag6.2
$$
for some $\ell ,m\in\bold N$.  By this method we produce from our
collection $\{z_n\}$ a new collection $\{z^*_n\}$.  Notice that it is
possible that $z^*_n=z^*_m$ even if $n\neq m$, but then $d(z_n)\sim
d(z_m)\sim\vert z_n-z_m\vert$.

\subheading{Step 2.  Another Geometric Construction}

Let $\{z_n\}$ be the collection of points in Step 1.  Let $\{ I_n\}$
be a collection of subarcs of
$\partial\tilde\Omega_k\backslash\partial\Omega_k$ such that
$\underset n\to\bigcup
I_n=\partial\tilde\Omega_k\backslash\partial\Omega_k$, $I_n\cap
I_m=\varphi$ when $n\neq m$, diameter$(I_n)\sim d(z_n)$, and $\vert
z-z_n\vert\leq\frac{1}{2} d(z_n)$ for $z\in I_n$.  We also define
$I_n^*$ to be the arc
$$
f^{-1}(\{(1-2^{-\ell})\text{exp}\{i(m+t)2^{-\ell}\pi\}:0\leq t\leq
1\}).
$$
See (6.2) for notation.  Then $I^*_n$ has diameter $\sim d(z_n)$ and
if $z\in I_n^*$ the hyperbolic distance from $z$ to $z_n^*$ (in
$\Omega_k$) is bounded by $C$.

With the notation of (6.1) we also denote by $J_n$ the subarc of
$\bold T$
$$
J_n=\{e^{i\theta}:m2^{-\ell}\pi <\theta\leq (m+1)2^{-\ell}\pi\}
$$
and we denote by $Q_n$ the ``square''
$$
Q_n=\{re^{i\theta}:(1-2^{-\ell})\leq r\leq 1, e^{i\theta}\in J_n\}.
$$

We now use the standard terminology that an arc $J_m$ is maximal in a
subcollection ${\Cal F}$ of $\{J_n\}$ if $J_m\in{\Cal F}$ and
$J_\ell\in{\Cal F}, \ell\neq m$, implies either $J_\ell\cap
J_m=\emptyset$ or $J_\ell\subset J_m$.  Notice (by the John condition)
that if $J_\ell\subset J_m$,
$$
\vert z_\ell -z^*_m\vert\leq Cd(z_m).\tag6.3
$$

Finally, if $\hat{\Cal F}$ is a subcollection of $\{I_n\}$ and
$E=\underset I_n\in\hat{\Cal F}\to\bigcup I_n$ we denote by $E^*$ the
set
$$
E^*=\underset J_n\in{\Cal F}\to\bigcup I^*_n
$$
where
$$
{\Cal F}=\{J_n:I_n\in\hat{\Cal F}~\text{and}~J_n~\text{is maximal}\}.
$$

\subheading{Step 3.  A Capacitary Estimate}

Let $E$ and $E^*$ be sets as in the previous paragraph.

\proclaim{Lemma 6.2}  $\text{Cap}(E^*)\geq c~\text{Cap}(E)$.
\endproclaim

\demo{Proof}  Let $\mu$ be a probability measure on $E$ satisfying
$$
\int\log\frac{1}{\vert z-\zeta\vert}d\mu (\zeta
)\leq\gamma~,~z\in\bold C.
$$
We relabel the intervals $J_n\in{\Cal F}$ so that ${\Cal
F}=\{J_1,J_2,\dots\}$ and $d(z_1)\geq d(z_2)\geq\cdots\geq d(z_n)\geq
d(z_{n+1}\geq\cdots$.  Define
$$
E_n=\{z\in E:\vert z-z_n\vert\leq Cd(z_n)~\text{and}~\vert z-z_m\vert
>Cd(z_m),m<n\},
$$
so that by (6.3),
$$
E=\underset n\to\bigcup E_n.
$$
Notice that the sets $E_n$ are pairwise disjoint.

Now define a probability measure $\mu^*$ by setting $\mu^*$ to have
uniform distribution on the center half $\frac{1}{2}I^*_n$ of $I^*_n,
\mu^*(I_n^*\backslash\frac{1}{2}I^*_n)=0$, and
$$
\mu^* (I^*_n)=\mu (E_n).
$$
By the construction of $I^*_n$ and $\frac{1}{2} I^*_n$,
$$
\text{dist}(\frac{1}{2} I^*_n,\frac{1}{2} I^*_m)\geq
cd(z_n)~,~\forall n\neq m.
$$
Let $z\in E$ and $z'\in\frac{1}{2} I_n$ where $n$ satisfies  $z^*\in
I_n$.  Then
$$
\align
\int\log\frac{1}{\vert z'-\zeta\vert}d\mu^*(\zeta ) &=\int_{\{\vert
z'-\zeta\vert\leq Ad(z)\}}+\int_{\{\vert z'-\zeta\vert >Ad(z)\}}\\
&\leq c+\int_{\{\vert z'-\zeta\vert\leq Ad(z)\}}\log\frac{1}{\vert
z-\zeta\vert}d\mu (\zeta )\\
&+c+\int_{\{\vert z'-\zeta\vert>Ad(z)\}}\log\frac{1}{\vert
z'-\zeta\vert}d\mu (\zeta )\\
&\leq 2c+\gamma,
\endalign
$$
and Lemma 6.2 is established.

\enddemo
\subheading{Step 4.  Proof of the Proposition}

By Lemmata 2.3 and 2.5 and by estimate (6.1),
$$
\vert H(z)-H(z^*)\vert\leq 1.
$$
Now let ${\Cal D} =\Omega_k\backslash\underset{}\Cal F
\to\bigcup\hat Q_n$, where the $Q_n$ are as defined in Step 2 and
$\hat Q_n=f^{-1}(Q_n)$.  Then by Lemma 4.1, ${\Cal D}$ is an
$(\varepsilon ')$ John domain.  We define our collection ${\Cal F}$
to be $\{I_n:\exists z\in I_n,\vert H(z)-H(z_k)\vert\geq\lambda\}$.
Then if $I_n\in{\Cal F},\vert H(\zeta )-H(z_k)\vert\geq\lambda -c$
for all $\zeta\in I_n$.  (This is why we slightly enlarges
$\tilde\Omega_k$ to $\Tilde{\Tilde\Omega}_k$.)  By our previous
estimate,
$$
\vert H(z)-H(z_k)\vert\geq\lambda -2c~\text{on}~I^*_n,
$$
for any $I_n\in{\Cal F}$.  Setting as before
$$
\align
E &=\underset{\Cal F}\to\bigcup I_n\\
\text{and}\qquad E^*&=\bigcup I_m^*,
\endalign
$$
we have $\text{Cap}(E^*)\geq c~\text{Cap}(E)$.  Now since
$$
\iint_{\Cal D}\vert\nabla H\vert^2dxdy\leq 1,
$$
it follows from Lemmata 2.3 and 2.5 that
$$
\text{Cap}(E^*)=o(1)~\text{as}~\lambda\to\infty.
$$
This completes the proof of Proposition 6.1.

\vskip .25in
\centerline{\bf \S7.  Proof of Theorem 1}

Let $\tilde\Omega_j$ and $\tilde\Omega_k$ be two domains satisfying
$$
\text{distance}(\tilde\Omega_j,\tilde\Omega_k)\leq 1.
$$

It is an exercise to find domains $\tilde\Omega_{j_{1}},\dots
,\tilde\Omega_{j_{N}}$ where $j_1=j,j_N=k$, and
$\partial\tilde\Omega_{j_{m}}\cap\partial\tilde\Omega_{j_{m+1}}$ is
an arc of diameter $\geq c$, and $N\leq C$.  (Use the fact that
$\Omega$ is a John domain and each $\tilde\Omega_j$ is an $\eta$
quasicircle, i.e. Lemma 5.2.)  Then by Proposition 6.1,
$$
\align
\vert F(z_j)-F(z_k)\vert &\leq\sum_{m=1}^{N-1}\vert
F(z_{j_{m}})-F(z_{j_{m+1}})\vert\\
&\leq C\sum_{m=1}^N\left(\iint_{\Tilde{\Tilde\Omega}_{j_{m}\backslash
K}}\vert\nabla F\vert^2dxdy\right)^{1/2}\\
&\leq C'\left(\iint_{\{z\in K^{c}:\vert z-z_{j}\vert\leq
C\}}\vert\nabla F\vert^2dxdy\right)^{1/2}
\endalign
$$

We also notice by (6.1) that if $z\in\tilde\Omega_k$ and $d(z)\geq
1$,
$$
\vert F(z)-F(z_k)\vert\leq
C\left(\iint_{\Tilde{\Tilde\Omega}_{k}\backslash K}\vert\nabla
F\vert^2dxdy\right)^{1/2}.
$$
Putting our last two estimates together we see there is $\tilde F\in
C^\infty (\bold R^2)$ such that $\tilde F(z)=F(z)$ when $d(z)\geq 1$ and
$$
\iint_{\{z\in K^{c}:\vert d(z)\vert\leq 1\}}\vert\nabla\tilde
F\vert^2dxdy\leq C\iint_{\{z\in K^{c}:\vert d(z)\vert\leq
C\}}\vert\nabla F\vert^2dxdy.
$$
Here we are using the fact that, by the construction of the
$\tilde\Omega_k$, $\{z\in K^c:d(z)\leq 1\}\subset\underset k\to\bigcup\tilde\Omega_k$.
Since the John condition is dilation invariant, we may now build a
sequence $\tilde F_n\in C^\infty (\bold R^2)$ with $\tilde
F_n(z)=F(z)$ when $d(z)\geq\frac{1}{n}$ and
$$
\iint_{\{z\in K^{c}:d(z)\leq\frac{1}{n}\}}\vert\nabla\tilde
F_n\vert^2dxdy\leq C\iint_{\{z\in
K^{c}:d(z)\leq\frac{c}{n}\}}\vert\nabla F\vert^2dxdy.
$$
Since $\vert K\vert =0$, it follows that $F\in W^{1,2}(\bold R^2)$ and
$$
\iint_{\bold R^{2}}\vert\nabla F\vert^2dxdy=\iint_{K^{c}}\vert\nabla
F\vert^2dxdy.
$$

\vskip .25in

\centerline{\bf \S8.  Proof of Theorem 2}

Let $\Omega$ be an $(\varepsilon )$ John domain with compact boundary
$K$ of diameter one, let $\{Q_j\}$ denote the Whitney decomposition
of $\Omega$ into dyadic squares \cite{11}, and let $z_j$
 be the center
of $Q_j$.  Let $A=A(\varepsilon )$ be a large constant and define
${\Cal F}_n=\{z_j:A^{-n}\leq d(z_j)\leq A\}$.  It is an exercise with
the John condition to construct a connected graph $G_0$ such
that every edge in $G_0$ is of the form $[z_j,z_k]$ where $\partial
Q_j\cap\partial Q_k\neq\phi$, and where the vertices $V_0$ of
$G_0$ satisfy
$$
{\Cal F}_0\subset V_0\subset\{z_j:1\leq d(z_j)\leq A^2\}.
$$
We also build $G_0$ so that
\roster
\item"{(8.1)}"  If $d(z_j)\geq 1$ and $z_j\notin V_0$ then $Q_j$ is
in the unbounded component of $\bold C\backslash\underset z_k\in{\Cal
F}_0\to\bigcup Q_k$.
\endroster

It is now an exercise (with induction) to construct connected graphs
$G_n$ with the following properties:
\roster
\item"{(8.2)}"  Every edge in $G_n$ is of the form $[z_j,z_k]$ for some
$z_j,z_k\in{\Cal F}_{n+1}\cup V_0$ where $\partial Q_j\cap\partial
Q_k\neq\emptyset$.
\item"{(8.3)}"  Every $z_j\in{\Cal F}_n$ is in $V_n$, the vertices of
$G_n$.
\item"{(8.4)}"  If $\rho (z_j,z_k)$ is the graph distance on $G_n$,
$$
\inf_{z_{j}\in{\Cal F}_{n}}\rho (z_j,z_k)\leq C~,~z_k\in G_n.
$$
\item"{(8.5)}"  $V_n\subset V_{n+1}$
\endroster

Notice that we have chosen $G_0$ to be connected.  Let $z_0\in V_0$
be an extreme point of the (planar set) convex hull $(G_0)$.  We may
assume by induction that each $G_n$ is actually a {\it directed}
graph in the following sense.  Each edge $[z_j,z_k]$ is directed in
the sense that (perhaps switching $j$ and $k$)
$$
\rho (z_j,z_0)=\rho (z_k,z_0)+1.\tag8.6
$$
Such an edge is an outgoing edge from $z_j$.   It is not hard to see
that we may choose the $G_n$ so that
\roster
\item"{(8.7)}"  Each $z_j\neq z_0$ has exactly one outgoing edge.
\endroster

\proclaim{Lemma 8.1}  The graph $G_n$ is simply connected, i.e. it
contains no loops.
\endproclaim

\demo{Proof}  Suppose to the contrary that there is a loop in $G_n$.
Let $z_j$ be a vertex in the loop maximizing $\rho (z_j,z_0)$.  Then
$z_j$ has two outgoing edges (by (8.6)) and this contradicts
(8.7).\qed

Let $G=\underset n\to\lim G_n$ be the limiting graph, so that $G$ is
simply connected.  It is clear that $K\cup G$ is connected.  Notice
by (8.3) that
\roster
\item{(8.7)}  For every $z_j\in G$ there is an arc $\gamma\subset G$
from $z_j$ to $z_0$ which satisfies the $\varepsilon '$ John
condition in $\Omega$.
\endroster
\flushpar In other words, $G$ is a John graph.

For a Whitney square $Q_j$ with $z_j\in G$ let $\{{\Cal L}_k^j\}$
denote all the edges of $G$ with one endpoint being $z_j$.  Define
$$
I_k^j=\{z\in\partial Q_j:~\text{distance}(z,{\Cal
L}_k^j)<\delta~\text{diam}(Q_j)\},
$$
where $\delta$ is a small constant, and put
$$
S_j=\partial Q_j\backslash\underset k\to\bigcup I_k^j~,~j\neq 0.
$$

For the special point $z_0\in G$ we select a Whitney square $Q_\ell$
such that $z_\ell\notin G$, $\partial Q_0\cap\partial
Q_\ell\neq\phi$, and we put
$$
\align
S_0 &=\partial Q_0\backslash (I_\ell^0\cup\underset k\to\bigcup
I_k^0),\\
\hat\Omega &=\Omega\backslash\underset z_j\in G\to\bigcup S_j.
\endalign
$$

\proclaim{Lemma 8.2}  $\hat\Omega$ is simply connected.
\endproclaim

\demo{Proof}  Let $\hat\Omega_+=\cup\{\hat\Omega\cap Q_j:z_j\notin
G\}$, $\hat\Omega_-=\cup\{\hat\Omega\cap Q_j:z_j\in G\}$ so that
$$
\hat\Omega =\hat\Omega_+\cup\hat\Omega_-\cup I_\ell^0.
$$
By condition (8.1), $\hat\Omega_+$ is simply connected (in $\bar{\bold
C}$), so it is only necessary to check that $\bar\Omega_-$ is simply
connected.

We first verify that $\hat\Omega_-$ is connected.  Let $z\in
Q_j\cap\hat\Omega_-$ and let $\gamma$ be an arc in $G$ connecting
$z_j$ to $z_0$.  Then $\gamma '=[z,z_j]\cup\gamma$ is an arc in
$\hat\Omega_-$ which connects $z$ to $z_0$.

Now suppose that $\gamma$ is a loop in $\hat\Omega_-$ that is not
homologous to zero.  It is then an elementary exercise to homotopy
$\gamma$ to $\gamma '$, a loop in $G$ that is not homologous to zero.
 This contradicts Lemma 8.1.\qed

It is clear from the construction of $\hat\Omega$ that
$\partial\Omega\subset\partial\hat\Omega$.  To verify that
$\hat\Omega$ is a John domain we must look at two cases.

\enddemo
\enddemo
\subheading{Case 1}  $z\in\hat\Omega_+\cap Q_j$.  There is arc
$\gamma$ from $z_j$ to some $z_k\notin G$ such that length$(\gamma
)\leq C$, $d(z)\geq 1$ on $\gamma$, and $z_k$ is not in the convex
hull of $\partial\hat\Omega$.  By selecting a suitable ray $R$ from
$z_k$ to $\infty$ we then see that
$$
[z,z_j]\cup\gamma\cup R
$$
is the required John arc.

\subheading{Case 2}  $z\in\hat\Omega_-\cup I_\ell^0$.  Let $\gamma$
be a John arc from $z_\ell$ (the center of the special Whitney
square $Q_\ell$ adjacent to $z_0$) to $\infty$.  Then if $z\in Q_j$
and $\gamma '\subset G$ is the John arc from $z_j$ to $z_0$ guaranteed by
condition (8.7), we see that
$$
[z,z_j]\cup \gamma '\cup [z_0,z_\ell]\cup\gamma
$$
is the required John arc.

\end